\documentclass[]{article}
\usepackage{latexsym,amssymb}
\title{Carter Subgroups, Amalgams, Simple Groups, and the $Z_{p}^{\ast}$-theorem}
\author{Geoffrey R. Robinson,\\ Institute of Mathematics,\\ University of Aberdeen,\\ Aberdeen AB24 3UE}

\begin{document}

\maketitle

\begin{abstract}
We consider an amalgam of groups constructed from fusion systems for different odd primes $p$ and $q$.
This amalgam contains a self-normalizing cyclic subgroup of order $pq$ and isolated elements of order $p$ and $q$.
\end{abstract}

\section{Introduction}

In earlier work, ([3],[4]) we used an (iterated) amalgam $X = X_{\mathcal{F}}$ to realise an Alperin fusion system $\mathcal{F}$ on a finite $p$-group $P$ via conjugation within $X$, and to obtain explicit  linear representations of $X$, thereby relating many finite groups (finite homomorphic images of $X$) to the original fusion system $\mathcal{F}$. 

\medskip
To partly motivate what follows, we outline an extension of the example in [4] which will also illustrate 
some ideas behind the main construction here. If we consider the maximal fusion system on a semidihedral 
$2$-group of order $16,$ this fusion system is realised by an amalgam $X = {\rm GL}(2,3)*_{D}S_{4}$ 
where $D$ is a dihedral group of order $8.$ It is interesting to note that $X$ has a unique conjugacy 
class of involutions, and that we have $C_{X}(t) = {\rm GL}(2,3)$ for each involution $t \in X$. There 
are only two non-isomorphic finite simple groups $G$ which contain an involution $u$ with $C_{G}(u) \cong 
{\rm GL}(2,3)$ (this is a theorem of R. Brauer). 

\medskip
The amalgam $X$ is certainly not a simple group, but all its proper non-trivial normal subgroups 
are free (as is implicitly noted in [4]). In [4], we showed that for each odd prime $q$, there is 
an epimorphism from $X$ to ${\rm SL}(3,q)$ (when $q \equiv 1,3$ (mod $8$)) or ${\rm SU}(3,q)$ 
(when $q \equiv 5, 7$ (mod 8)). Hence this fusion system on this single $2$-group leads naturally 
to an infinite number of non-isomorphic finite simple groups (not all of which have a semi-dihedral 
Sylow $2$-subgroup of order $16$). Furthermore, each epimorphism constructed has free kernel, so 
the finite groups used to build the amalgam embed faithfully in each of these simple epimorphic 
images, and these epimorphic images are all generated by the images of these finite groups. 

\medskip 
In this paper, then, we will construct an amalgam $X$ from two different fusion systems, for 
two different primes $p$ and $q.$ Again, the construction illustrates a general methodology 
which should be much more widely applicable (and with iterated amalgams). In the case we consider 
below, each fusion system is of a rather transparent type (and is a fusion system on an extra-special 
group) and the interaction between the fusion systems is almost minimal. 

\medskip
However, it illustrates the general methodology, and (in our view) it also illustrates that rich 
structures can arise in this context from uncomplicated building blocks. We were initially led 
to these considerations by a realisation of a connection between odd analogues of Glauberman's 
$Z^{\ast}$-theorem, [2], and some troublesome configurations considered in the Feit-Thompson proof 
of the solvability of finite groups of odd order, [1]. 

\medskip
We will construct amalgams $X$ realising related configurations within perfect infinite 
groups, and show that all proper normal subgroups of the amalgams constructed are free. 
Each amalgam $X$ constructed is a perfect group which realises a constrained fusion 
system on an extra-special (finite) $p$-group $P$ of exponent $p$ and a constrained fusion system 
on an extra-special (finite) $q$-group $Q$ of exponent $q$ where $p$ and $q$ are distinct odd primes. 
Furthermore, each element of $Z(P)$ and each element of $Z(Q)$ is isolated in $X.$ The amalgams 
$X$ all have free normal subgroups of finite index, so do have non-Abelian finite simple groups 
as epimorphic images, which are generated by the images of $N_{X}(P)$ and $N_{X}(Q)$. We 
exhibit an explicit pair of generators for each amalgam $X$, one of order $pq$ and one of infinite 
order. Also, each $X$ is generated by three explicitly identified elements, of respective orders 
$q,p$ and $pq.$ Hence we have explicit generators for any finite simple homomorphic image of $X,$ 
though the order of the image of the element of infinite order is not a priori obvious. 

\section{Notation, Definitions, Background} 

A finite $p$-subgroup $S$ of a (possibly infinite) group $G$ is said to be a Sylow $p$-subgroup 
of $G$ if every finite $p$-subgroup of $G$ is conjugate to a subgroup of $S.$ We recall that 
a Carter subgroup of a finite solvable group $H$ is a self-normalizing nilpotent subgroup of $H.$ 
The fact that such an $H$ always has a Carter subgroup, and that these are all $H$-conjugate, 
was proved by R. Carter. 

\medskip
We refer to an element $x$ of a group $G$ as \emph{isolated} if $x$ commutes with none of its other 
$G$-conjugates. This terminology occurs frequently in existing literature, though usually in 
the context of finite groups. Glauberman's $Z^{\ast}$-theorem proves that if $G$ is a finite 
group with no non-trivial normal subgroup of odd order, then any isolated involution $t \in G$
is central in $G$.

\medskip
 For ease of later discussion, we will call a group $H$ \emph{torsion simple} if all its proper 
non-trivial normal subgroups are torsion free. In the amalgams we deal with in this paper, 
torsion free subgroups are always free. 

\medskip
In the Feit-Thompson proof of the solvability of finite groups of odd order,[1], it is shown that 
a putative minimal finite simple group of odd order must have a self-normalizing cyclic subgroup 
of order $pq$ for distinct (odd) primes 
$p$ and $q$ (and that there is a unique conjugacy class of such self-normalizing cyclic subgroups). 
This self-normalizing cyclic subgroup is the intersection of a pair of maximal subgroups. 
One of these has an elementary Abelian normal Sylow $p$-group which is the Frobenius kernel of 
a Frobenius subgroup of index $q$, and the other has an elementary Abelian normal Sylow $q$-subgroup which is the Frobenius kernel of a Frobenius subgroup of index $p.$ 

\medskip
In the Odd Order paper,[1], this possibility is eventually eliminated by a difficult analysis 
in Chapter VI (the last chapter).

\section{The Construction} 

Let $p$ and $q$ be distinct odd primes. Let $m$ be the smallest positive integer such that 
$|{\rm Sp}(2m, q)| $ is divisible by $p$ and let $n$ be the smallest positive integer such 
that $|{\rm Sp}(2n, p)|$ is divisible by $q$. Let $A$ be the semi-direct product  of an 
extra-special $q$-group $Q$ of order $q^{2m+1}$ and exponent $q$ with a cyclic group of 
order $p$ acting trivially on $Z(Q)$ and faithfully on $Q/Z(Q)$. 
Let $B$ be the semi-direct product  of an 
extra-special $p$-group $P$ of order $p^{2n+1}$ and exponent $p$ with a cyclic group of 
order $q$ acting trivially on $Z(P)$ and faithfully on $P/Z(P)$. 

\medskip
Notice that $A$ has a Carter subgroup which is cyclic of order $pq$ and contains 
$Z(Q)$, while $B$ has a Carter subgroup which is cyclic of order $pq$ and contains 
$Z(P).$ 

\medskip
We form the amalgam $X = A *_{C}B,$ where $C$ is a cyclic group of order $pq$ identified 
(by inclusion) with a Carter subgroup of $A$ and a Carter subgroup of $B.$ Then $Q$ is a 
Sylow $q$-subgroup of $X$ and $P$ is a Sylow $p$-subgroup of $X.$ 

\medskip
With the above identification, we take $C = Z(P) \times  Z(Q).$ Furthermore, as will become 
apparent below, the fusion system induced on $Q$ by conjugation within $X$ is just the constrained 
fusion system induced by $A$ on $Q$ and the fusion system induced on $P$ by conjugation within 
$X$ is just the constrained fusion system induced on $P$ by $B.$ This could be proved from the 
results of [3], but we will give a self-contained proof below. 

\medskip
We claim that $X$ is torsion simple. Let $N$ be a proper non-trivial normal subgroup of 
$X.$ If $N$ is not free, then (under current hypotheses) it is not torsion free either, and we 
either have $N \cap A \neq 1$ or $N \cap B \neq 1,$  since all elements of finite order in $X$ 
lie in a conjugate of $A$, or a conjugate of $B$.  Hence $N$ contains either an element of order 
$p$ or an element of order $q.$ Suppose that $N$ contains an element of order $p.$ If this element 
lies in $B,$ then $N \cap  Z(P) \neq 1$ and $[Q,Z(P)] = Q \leq N.$ Then also $[P,Z(Q)] = P \leq N,$ so $N = X.$ 
A similar argument gives $N = X$ if $N$ contains an element of order $q.$ Hence $N$ is free. 

\medskip 
Now we must have $X = [X,X],$ for otherwise $[X, X]$ is free, so that $A$ and $B$ each embed isomorphically into $X/[X,X]$, a contradiction, as both $A$ and $B$ are non-Abelian. As explained in [4], $X$ has a free 
(not necessarily normal) subgroup of index $|P ||Q|,$ so $X$ does have free normal subgroups of finite index. 
We note in passing that whenever $F$ is a maximal free normal subgroup of finite index in $X,$ then $X/F$ 
acts faithfully on $F/[F,F]$. This is a standard argument, but we reproduce it now: let $N =[F, X] \geq [F,F].$ 
Then $X/N$ is a central extension of $X/F$ by the finitely generated Abelian group $F/N.$ Then $F/N$  must 
be finite, otherwise there is a perfect central extension of the finite group $X/F$ by  $Z/rZ$ for every prime $r, $ a contradiction. In particular, $N> [F,F],$ since $F/[F,F]$ is infinite. Hence $C_{X}(F/[F,F])$ is a proper normal subgroup of $X$, hence free. By the maximal choice of $F,$ we must have $F = C_{X}(F/[F,F])$.

\medskip
Now we note that $C$ is a self-normalizing cyclic subgroup of $X.$ It is a general fact about amalgams 
that every element of $X\backslash C$ may be written in the form $u_{1}u_{2} \ldots u_{n},$ where  
each $u_{i}$ is either in $A \backslash C$ or in $B\backslash C,$ and there is no value of $i$ such 
that both $u_{i}$ and $u_{i+1}$ lie in $A$, nor is there any value of $j$ such 
that both $u_{j}$ and $u_{j+1}$ lie in $B$. Furthermore, each such product does lie outside $C.$ 
Now let $d$ be a generator of $C.$ Consider $u_{n}^{-1} \ldots u_{1}^{-1} d u_{1}\ldots u_{n}$ and set $v_{i} = d^{-1}u_{i} d$ for each $i.$ Then $u_{n}^{-1} \ldots u_{1}^{-1} d u_{1}\ldots u_{n} = dv_{n}^{-1} \ldots (v_{1}^{-1} u_{1}) \ldots u_{n}.$ Now $u_{i}$ does not normalize $C,$ since $C$ is self-normalizing in both 
$A$ and $B.$ In particular, $v_{1}^{-1}u_{1} \not \in C,$ but does lie in the same member of $\{A,B \}$ as $u_{1}$ does. Hence the given conjugate of $d$ lies outside $C,$ since $v_{j}$ lies outside $C$ for 
$j> 1,$ but does lie in the same member of $\{A,B\}$ as $u_{j}$ does. 

\medskip
Now we claim that $A \cap A^{x}$ is a $q^{\prime}$-group for each $x \in X \backslash A$
(and similarly $B \cap B^{y}$  is a $p^{\prime}$-group for each $y \in X \backslash B).$
Consider such an element of the form $x = u_{1}u_{2} \ldots u_{n}$ with each $u_{i}$ lying 
outside $C$ but inside $A$ or $B,$ such that there is no value of $i$ for which both 
$u_{i}$ and $u_{i+1}$ lie in the same member of $\{A,B\}$.  If possible, choose an element 
$a$ of order $q$ in $A \cap A^{x}.$ If $n =1,$ then $x = u_{1} \in B \backslash C,$ so $x^{-1}ax$ 
lies outside $A \cup B$ unless $a \in C$. But if $a \in C,$ then $\langle a \rangle = Z(Q),$ and $x^{-1}ax$
can only lie in  $A$ if it lies in $A \cap B = C.$ In that case, $x$ must lie in $N_{B}(Z(Q)) = C$, a contradiction.

\medskip
We note that if $n> 1,$ then $x$ lies outside $A \cup B$. Suppose now that $n>1 $ and consider 
$u_{n}^{-1} \ldots u_{1}^{-1} a u_{1} \ldots u_{n}.$ If $a \not \in C,$ this product lies outside 
$A \cup B$, for if $u_{1} \in A,$ then we may bracket it as\\  
$$u_{n}^{-1} \ldots (u_{1}^{-1}au_{1}) \ldots u_{n}.$$ 
Note that $u_{1}^{-1}au_{1}$ still lies in $A \backslash C.$ On the other hand, if $u_{1} \in B,$  the expression\\
$u_{n}^{-1} \ldots u_{1}^{-1} a u_{1} \ldots u_{n}$ is already expressed as a product of elements 
which lie alternately in $A \backslash C$ or $B \backslash C.$ 

\medskip
However, if  $a \in C,$ then we might as well suppose that $u_{1} \in  B \backslash C,$ since $O_{q}(C)= Z(A).$ 
But then, as before, we can express the product as $$a[a^{-1}u_{n}^{-1}a] \ldots [a^{-1}u_{1}^{-1} a u_{1}]u_{2} \ldots u_{n}.$$ 
Now notice that $[a, u_{1}] \in B.$ If $[a, u_{1}] \in C,$ then as above, we obtain $u_{1} \in  
N_{B}(Z(Q)) = C,$ contrary to assumption. Thus $[a, u_{1}] \in B \backslash C,$ which shows that the stated 
product lies outside $A \cup B.$ 

\medskip
Now we have proved that $A \cap A^{x}$  is a $q^{\prime}$-group for each $x \in X\backslash A$, and by 
symmetry, $B \cap B^{y}$ is a $p^{\prime}$-group for each $y \in X \backslash B.$ We now note that $z$ 
is isolated whenever $z \in Z(Q)^{\#}$,  (and an analogous statement holds for any $w \in Z(P )^{\#}).$ 
For suppose that $z^{x} \neq z$ commutes with $z.$ Then $\langle z, z^{x} \rangle$ is finite of order $q^{2}$
so is conjugate to a subgroup of $Q.$ Then $\langle z^{y},z^{xy} \rangle \leq Q$ for some $y \in X.$ Then 
$y \in A$ since $z \in A \cap yAy^{-1}$.  Similarly $xy \in A,$ so that $x \in A = C_{X}(z)$ and $z = z^{x}.$ Notice that if $N$ is any proper normal subgroup of $X,$ then $z \not \in N$ as $N$ is free, but $zN$ is not central in $X/N.$ For $X/N$ is not Abelian as $X = [X,X]$, so if $zN$ is central in $X/N,$ then $z$ lies in a proper normal subgroup of $X,$ which is free, a contradiction.

\medskip
\noindent {\bf Remark}: We note that $X$ (and hence every homomorphic image of $X$) is generated by two 
elements. For any $u \in Q \backslash Z(Q),$ and $v \in P \backslash Z(P),$ we have $X = Y = \langle uv, C \rangle.$ For let $z$ be a generator of $Z(Q)$ and $w$ be a generator of $Z(P ).$ Then $uv^{z} \in Y,$ so $v^{-1}v^{z} = (uv)^{-1}(uv)^{z} \in Y.$ Now $v^{-1}v^{z} \in P \backslash Z(P),$ so that $P =\langle v^{-1}v^{z},z \rangle$ and $P \leq Y.$ Similarly, $Q \leq Y,$ so $Y = X.$ Notice, however, that $uv$ 
has infinite order. We can certainly generate $X$ by three elements of finite order, as $X = \langle
u, v, zw \rangle$. 

\medskip
In fact, we claim that $T = \langle z^{b}w, uv\rangle = X$ for any $b \in B.$ For 
$T$ contains $\langle w, uv \rangle = A,$ from the above argument. Hence $T$ contains
$\langle z^{b},v \rangle$ which is a $B$-conjugate of $\langle z, bvb^{-1}\rangle = B$
as above. Thus $T = X.$ 

\medskip
We also remark that every element lying outside all conjugates of $A$ and outside all conjugates 
of $B$ has a centralizer which is free of rank $1.$ For if $x \in X$ and no conjugate of 
$x$ lies in $A \cup B$ then $x$ has infinite order. Hence $x \not \in C_{X}(a)$ for any $a \in A^{\#},$ 
and similarly for any conjugate of $A.$ A similar argument holds for conjugates of $B.$ Thus $C_{X}(x)$ 
is torsion free, and hence free. But a free group with non-trivial centre is free of rank $1.$ 

\medskip
We note that (by essentially the same proof as given in [4]), all normal subgroups of finite index 
of $X$ are free of rank greater than $1.$ For let $N$ be a normal subgroup of $X$ of finite 
index. We know already that $N$ is free, so that $A$ and $B$ embed isomorphically into $X/N$. Hence $[X:N]$ id divisible by $p^{2n+1}q^{2m+1}.$ If $N$ has $r$ generators, then $\chi(N) = \frac{1-r}{[X:N]}.$
We also have $$\chi(X) = (\frac{1}{qp^{2n+1}}+\frac{1}{pq^{2m+1}} - \frac{1}{pq}),$$ which yields 
$$r = 1 + \frac{(p^{2n}q^{2m} - p^{2n}-q^{2m})[X:N]}{p^{2n}q^{2m}} > 1.$$

\medskip
In particular, this implies that every non-identity element of $X$ has an infinite number of conjugates.
Hence over any field, the group algebra $FX$ ( consisting of finite $F$-linear combinations of elements of $X$)
has one-dimensional centre.

\section{Remarks on the finite case}

Although much more sophisticated counting arguments and character-theoretic arguments are used in Chapter 
V of ([1]), to reach the configuration which is left to be eliminated in Chapter VI of that work,  as a 
matter of interest, we give a direct proof here of an easier result. We prove that the $p$-local and $q$-local structure which we have shown to be present in $X$ can't all be present in a finite group.

\medskip
\noindent {\bf Theorem :} \emph{There is no finite group $G$ with the following properties:} 
\medskip
\noindent \emph{i) There are distinct odd prime divisors $p$ and $q$ of $|G|$ such that 
$G$ has a self-normalizing cyclic subgroup $C$ of order $pq$ with $C = Z(P) \times Z(Q)$ 
for $P$ a Sylow $p$-subgroup and $Q$ a Sylow $q$-subgroup of $G.$} 

\medskip
\noindent \emph{ii) $N_{G}(P) =  PZ(Q)= C_{G}(Z(P ))$ and  $N_{G}(Q) =  QZ(P)= C_{G}(Z(Q )).$ }
 
\medskip
\noindent \emph{iii) $P \cap P^{g} = 1 $ for all $g \in G \backslash N_{G}(P)$ and  $Q \cap Q^{h} = 1 $ for all $h \in G \backslash N_{G}(Q).$}

\medskip
\noindent {\bf Proof: } Suppose otherwise, and set $A = N_{G}(P)$, $B = N_{G}(Q)$. 
We may, and do, suppose that $G = \langle P, Q \rangle .$ Let $N$ be a normal subgroup of $G.$ 
Then $A \cap N \neq 1$ or $B \cap N = 1$ gives $Z(P) \leq N$ or $Z(Q) \leq N$, so either $Q = [Q,Z(P)] \leq N$ or $P = [P,Z(Q)] \leq N.$ Suppose that that $N$ is a $\{p,q\}^{\prime}$-group. Now $P$ is not cyclic, otherwise $G$ has a normal $p$-complement, whereas $P = [P,Z(Q)] \leq G^{\prime}$, and likewise $Q$ is not cyclic.
Now  $N\leq \langle C_{G}(x): x \in P^{\#}  \rangle \leq A,$ so that $N = 1$. Hence $G$ is a simple group. 
Now $G$ has $[G : A](|P |- 1)$ non-identity $p$-elements and $[G : B](|Q|- 1)$ non-identity $q$-elements. 
Also $G$ has (at least) $|G|\frac{(p-1)(q-1)}{pq}$ elements of order $pq$. This accounts for 
$$|G|[ \frac{1}{p} + \frac{1}{q} - \frac{1}{q|P|} - \frac{1}{p|Q|} + 1 -\frac{1}{p} - \frac{1}{q} + \frac{1}{pq}] $$ elements. However, this is greater than $|G|$ as $p$ and $q$ are both odd, a contradiction.

\begin{center}
{\bf Bibliography}
\end{center}

\medskip
\noindent [1] Feit, Walter; Thompson, John G., \emph{ Solvability of groups of odd order}, 
Pacific J. Math., {\bf 13}, (1963), 775-1029. 

\medskip
\noindent [2] Glauberman, George, \emph{Central elements in core-free groups}, J. Algebra, {\bf 4}, (1966), 403-420. 

\medskip
\noindent [3] Robinson, Geoffrey R., \emph{Amalgams, blocks, weights, fusion systems and finite simple groups}, J. Algebra, 314,{\bf 2}, (2007), 912-923.

\medskip
\noindent [4] Robinson, Geoffrey R., \emph{Reduction mod q of fusion system amalgams}, Trans. Amer. Math. Soc. 363, {\bf 2}, (2011), 1023-1040.

\end{document}